\documentclass[11pt,twoside,a4paper]{amsart}[2000]
\title{Rokhlin actions and self-absorbing $C^*$-algebras}
\author{Ilan Hirshberg}
\address{Department of Mathematics, Ben Gurion University of the Negev, P.O.B. 653, Be'er \indent Sheva 84105, Israel}
\email{ilan@imada.sdu.dk}

\author{Wilhelm Winter}
\address{Mathematisches Institut der Universit\"at M\"unster\\
Einsteinstr. 62\\ 48149 M\"unster,
\indent Germany}

\thanks{The  second
author was supported  by the DFG (through the SFB 478) and the EU-Network  \indent Quantum  Spaces - Noncommutative
Geometry (Contract No. HPRN-CT-2002-00280).}

\email{wwinter@math.uni-muenster.de}

\usepackage{amssymb}
\usepackage{amsthm}
\usepackage{xypic}

\theoremstyle{plain}

\newtheorem{Thm}{Theorem}[section]
\newtheorem{Cor}[Thm]{Corollary}
\newtheorem{Lemma}[Thm]{Lemma}

\theoremstyle{definition}
\newtheorem{Def}[Thm]{Definition}
\newtheorem{Not}[Thm]{Notation}

\newcommand{\B}{\mathcal{B}}
\newcommand{\A}{\mathcal{A}}
\newcommand{\M}{\mathcal{M}}

\newcommand{\D}{\mathcal{D}}
\newcommand{\Ch}{\mathcal{C}}
\newcommand{\U}{\mathcal{U}}
\newcommand{\Zh}{\mathcal{Z}}

\newcommand{\T}{{\mathbb T}}
\newcommand{\R}{{\mathbb R}}
\newcommand{\N}{{\mathbb N}}
\newcommand{\Z}{{\mathbb Z}}
\newcommand{\C}{{\mathbb C}}

\renewcommand{\arrow}{\rightarrow}

\newcommand{\aut}{\mathrm{Aut}}

\newcommand{\id}{\mathrm{id}}

\newcommand{\OO}{\mathcal{O}}

\newcommand{\eps}{\varepsilon}

\newcommand{\Ainf}{\A_{\infty}}
\newcommand{\Acent}{\A_{\infty}\cap \A'}
\newcommand{\MAcent}{\M(\A)_{\infty}\cap \A'}

\newcommand{\Ainfc}{\A_{\infty}^{(\alpha)}}
\newcommand{\lNAc}{\ell^{\infty,(\alpha)}(\N,\A)}
\newcommand{\Acentc}{\Ainfc\cap \A'}

\newcommand{\oinn}{\overline{\mathrm{Inn}_0}}

\date{\today}


\begin{document}

\begin{abstract}
Let $\A$ be a unital separable $C^*$-algebra, and $\D$ a $K_1$-injective strongly
self-absorbing $C^*$-algebra. We
show that if $\A$ is $\D$-absorbing, then the crossed product of
$\A$ by a compact second countable group or by $\Z$ or by $\R$ is $\D$-absorbing as well, assuming the
action satisfying a Rokhlin property. In the case of a compact Rokhlin action we prove a similar statement about approximate divisibility.
\end{abstract}
\maketitle

\section{Introduction}

Following the terminology in a recent paper of A. Toms and the
second author (\cite{toms-winter}), we call a separable, unital
$C^*$-algebra $\D$ \emph{strongly self-absorbing} if it is
infinite-dimensional and the map $\D \arrow \D \otimes \D$ given
by $d \mapsto d \otimes 1$ is approximately unitarily equivalent
to an isomorphism (we note that strongly self-absorbing algebras
are always nuclear, so there is no ambiguity when we talk of
tensor products).
 Currently, the only known examples of such
algebras are the Jiang-Su algebra $\mathcal{Z}$ (\cite{jiang-su}), the Cuntz algebras $\OO_2$ and $\OO_{\infty}$, 
UHF algebras of infinite type (i.e. where all the primes which
occur in the relevant supernatural number do so with infinite
multiplicity)  and tensor products of
$\OO_{\infty}$ by such UHF-algebras. Those algebras exhaust the
possible Elliott invariants for strongly self-absorbing algebras,
and thus one might hope that this list is complete.

While only few algebras are strongly self-absorbing, many
$C^*$-algebras $\A$ are $\D$-absorbing for a strongly
self-absorbing algebra $\D$ (i.e. $\A \cong \A \otimes \D$), and
such algebras seem to enjoy nice regularity properties -- see
\cite{rordam UHF 1,rordam UHF 2} for absorption of UHF algebras,
and \cite{gong-jiang-su,rordam Z,toms-winter ASH} for absorption of the Jiang-Su
algebra. Absorption of $\OO_{\infty}$ and $\OO_2$ plays a central
role in the classification theorems of Kirchberg and Phillips
(\cite{kirchberg-phillips,kirchberg book,phillips 1}), and is the
focus of further study (see for instance \cite{kirchberg-rordam}).
It thus seems interesting to study the permanence properties of
$\D$-absorption. In \cite{toms-winter}, it was shown that if $\D$
is strongly self-absorbing and $K_1$-injective (i.e., the canonical map $\U(\D)/\U_{0}(\D) \to K_{1}(\D)$ is injective, a condition which is automatically fulfilled for the known examples mentioned above), then the property of being
(separable and) strongly self-absorbing is closed under passing to
hereditary subalgebras, quotients, inductive limits and extensions
(for $\D=\OO_{2}$ and $\D=\OO_{\infty}$ these results had already been shown by Kirchberg; see \cite{kirchberg}).

This note concerns the question of permanence under formation of
crossed products. One cannot expect permanence to hold in general.
Indeed, $\OO_2$ is $\OO_2$-absorbing, however there are actions of
$\Z_2 = \Z/2\Z$ on $\OO_2$ such that the crossed product has non-trivial
$K$-theory (see \cite{izumi}); in particular, the crossed product
algebra cannot be $\OO_2$-absorbing. Note that if $\alpha$ is such an
order 2 automorphism, then $\OO_2 \times_{\alpha} \Z_2$ is a
quotient of $\OO_2 \times_{\alpha} \Z$. Since $\OO_2$-absorption
passes to quotients, we see that $\OO_2$-absorption is not
permanent under crossed products by $\Z$ either.

In the present paper we show that this phenomenon does not occur if the group action satisfies certain extra conditions. More precisely, we show the following:

\begin{Thm}
\label{main result}
Let $\A$ and $\D$ be separable unital $C^{*}$-algebras; let $\D$ be strongly self-absorbing. Suppose $\alpha: G \to \aut (\A)$ is a strongly
continuous action of a group $G$, where $G$ is second countable
compact Hausdorff, $\Z$ or $\R$, which satisfies the respective
Rokhlin property (to be made precise in the subsequent sections).
\begin{enumerate}
\item If $G$ is compact and $\A$ is $\D$-absorbing, then  $\A \times_{\alpha}
G$ is $\D$-absorbing. If $\A$ is approximately divisible, then so is
$\A\times_{\alpha} G$.
\item If $G = \Z$ or $G = \R$, $\A$ is $\D$-absorbing and $\D$ is
$K_1$-injective, then $\A \times_{\alpha} G$ is $\D$-absorbing.
\end{enumerate}
\end{Thm}

There are many examples of $C^{*}$-algebras with Rokhlin actions;  see \cite{izumi,izumi II} for the finite group case and \cite{kishimoto, kishimoto UHF} for the case where $G=\Z$ or $\R$.  Rokhlin flows and a number of striking applications to the theory of purely infinite $C^{*}$-algebras have been studied in \cite{kishimoto flows} and \cite{BKR} and the references therein.    \cite{izumi survey} provides a  survey (mainly in the case where $G=\Z$ or $\R$) and describes a number of  applications to the theory of von Neumann algebras and to Elliott's classification program.  \\

The first author would like to thank Mikael R{\o}rdam for
suggesting this line of research, and for some helpful subsequent
conversations concerning this paper, as well as N. Christopher
Phillips for a helpful conversation. The second author is indebted
to Siegfried Echterhoff and Andrew Toms for a number of
enlightening comments and discussions.

\section{Central sequence embeddings and crossed products}

\begin{Not}
Let $\A$ be a $C^{*}$-algebra. We denote $$\A_{\infty} = \ell^{\infty}(\N,\A)/\Ch_0(\N,\A) \, .$$ $\A$ may be embedded  into   $\ell^{\infty}(\N,\A)$ and into $\A_{\infty}$ in a canonical way (as constant sequences); we shall write $\iota_{\A}$ for both these embeddings -- although we will sometimes find it convenient not to state them explicitly.\\
We write $\Acent$ for the central sequence algebra of $\A$, i.e., the relative commutant of $\A$ in $\A_{\infty}$.\\
If $\alpha:G \arrow \aut(\A)$ is a strongly continuous action of a
locally compact group $G$ on a $C^{*}$-algebra $\A$, then we have
naturally induced actions of $G$ on $\ell^{\infty}(\N,\A)$,
$\A_{\infty}$ and $\Acent$, respectively -- those actions will all
be denoted by $\bar{\alpha}$. They are in general not strongly
continuous. We denote by $\lNAc$ the set of elements of
$\ell^{\infty}(\N,\A)$ on which the action $\bar{\alpha}$ of $G$
is continuous -- this is clearly a $C^*$-algebra, which contains
$\A$ as the constant sequences, as well as $\Ch_0(\N,\A)$. We
denote $$\Ainfc = \lNAc / \Ch_0(\N,\A)$$ (we stress that $\Ainfc$
is not defined as the set of all elements in $\Ainf$ on which $G$
acts continuously -- it is a-priori smaller).
 We denote by $\M(\A)$ the multiplier algebra of
$\A$.  To avoid notational
nuisances, we follow in this paper the convention that $0$
is not considered to be a unital $C^*$-algebra (i.e., we require
that $1 \neq 0$).
\end{Not}

We have the following characterization of $\D$-absorption (based
on ideas of Elliott and of Kirchberg), which appears as Theorem
7.2.2 in \cite{rordam book}. We note that the statement in
\cite{rordam book} refers to the relative commutant of $\A$ in an
ultrapower of $\M(\A)$; however, it is easy to see that the
characterization still holds as stated below.

\begin{Thm}
\label{D-stable characterization} Let $\D$ be  a strongly
self-absorbing and $\A$ be any separable $C^{*}$-algebra. $\A$ is
$\D$-absorbing if and only if $\D$ admits a unital
$*$-homomorphism to $\MAcent$.
\end{Thm}
Note that since any strongly self-absorbing $C^*$-algebra has to
be simple, it follows that unless $\A = 0$, the $*$-homomorphism
above must be an embedding (of course, $0$ is $\D$-absorbing for
any $\D$).

The following gives us a sufficient condition for $\D$-absorption
of the crossed product of $\A$ by a group.

\begin{Lemma}
\label{invariant embedding implies absorption} Let $\A$, $\D$ be
unital separable $C^{*}$-algebras. Let $G$ be a locally compact
Hausdorff group with a strongly continuous action $\alpha: G
\arrow \aut(\A)$. Suppose there is an embedding $\D\arrow \Acent$
whose image is fixed under the induced action of $G$. Then, $\D$
admits a unital embedding into $(\M(\A \times_{\alpha}G))_{\infty}
\cap (\A \times_{\alpha}G)'$. In particular, if $\D$ is
strongly self-absorbing, then $\A \times_{\alpha} G$ is
$\D$-absorbing.
\end{Lemma}
\begin{proof}
By the universal property of $\A \times_{\alpha} G$, there are a unital 
$*$-homomorphism
\[
\pi:\A \to \M(\A \times_{\alpha} G)
\]
and a strictly continuous unitary representation
\[
u: G \to \M(\A \times_{\alpha}G)
\]
such that
\[
u_{g} \pi(a) u_{g^{-1}} = \pi \alpha_{g}(a) \; \forall \, a \in
\A, \, g \in G \, .
\]
$\pi$ and $u$ induce maps
\[
\tilde{\pi}: \ell^{\infty}(\N,\A) \to \ell^{\infty}(\N,\M(\A
\times_{\alpha} G)) \; , \; \tilde{u}: G \to
\ell^{\infty}(\N,\M(\A \times_{\alpha} G))
\]
which in turn induce maps
\[
\bar{\pi}:\A_{\infty} \to (\M(\A \times_{\alpha} G))_{\infty} \; ,
\;  \bar{u}: G \to (\M(\A \times_{\alpha}G))_{\infty}
\]

Let $\varphi:\D \arrow \Acent \subseteq \Ainf$ be the given
embedding, and consider $\bar{\pi} \circ \varphi : \D \arrow
(\M(\A \times_{\alpha} G))_{\infty}$. We claim that the image of the unital $*$-homomorphism 
$\bar{\pi} \circ \varphi$ commutes with the images of
$\bar{\pi}\circ\iota_{\A}$ and $\bar{u}$. That it commutes with
the image of $\bar{\pi}\circ\iota_{\A}$ is immediate. As for
$\bar{u}$, fix $g \in G$, $d \in \D$. Lift $\varphi(d)$ to
$(d_1,d_2,...) \in \ell^{\infty}(\N,\A)$. We have $u_g \pi(d_n)
u_{g^{-1}} = \pi(\alpha_g(d_n))$ for all $n$, so,
$\tilde{u}_g\tilde{\pi}(d_1,d_2,...)\tilde{u}_{g^{-1}} =
\tilde{\pi}(\alpha_g(d_1),\alpha_g(d_2),...)$. We know that
$$
(\alpha_g(d_1),\alpha_g(d_2),...) - (d_1,d_2,...) \in\Ch_0(\N,\A)
\; ,
$$
so
$$
\tilde{u}_g\tilde{\pi}(d_1,d_2,...)\tilde{u}_{g^{-1}} -
\tilde{\pi}(d_1,d_2,...) \in \Ch_0(\N,\M(\A \times_{\alpha} G))
$$
and thus $\bar{u}_g\bar{\pi}(\varphi(d))\bar{u}_{g^{-1}} -
\bar{\pi}(\varphi(d)) = 0$, as required.

Therefore, $\bar{\pi} \circ \varphi(\D)$ commutes with
the image of $\A \times_{\alpha} G$ in $(\M(\A \times_{\alpha}
G))_{\infty}$, as required. The second statement follows from
Theorem \ref{D-stable characterization}.
\end{proof}

In fact, to deduce $\D$-absorption, it is enough to show that the
conditions of Lemma \ref{invariant embedding implies absorption}
hold only approximately, as made precise in the following lemma.

\begin{Lemma}
\label{approximately fixed commuting embedding} Let $\A,\B$ be
unital separable $C^{*}$-algebras. Let $G$ be a second-countable
locally compact Hausdorff group with a strongly continuous action
$\alpha: G \arrow \aut(\A)$. Suppose that for any finite sets $B_0
\subseteq \B$, $A_0 \subseteq \A$, any compact subset $K_0
\subseteq G$ and any $\eps>0$ there is a completely positive
contraction $\varphi:\B\arrow \Ainfc$ such that for all $b,b' \in
B_0$, $a \in A_0$, $g \in K_0$ we have
\begin{itemize}
    \item[(i)] $\|\bar{\alpha}_g(\varphi(b)) - \varphi(b)\|<\eps$
    \item[(ii)] $\|\varphi(1)-1\|<\eps$
    \item[(iii)] $\|\varphi(b)\varphi(b') - \varphi(bb')\|<\eps$
    \item[(iv)] $\|[\varphi(b),a]\|<\eps$
    \item[(v)] $|\|b\|-\|\varphi(b)\| | <\eps$.
\end{itemize}
Then, there is an embedding $\B\arrow \Acent$ whose image is fixed
under the induced action of $G$. If $\B$ is simple, the conclusion holds without assuming condition (v) above.
\end{Lemma}
\begin{proof}
Pick increasing sequences of finite sets $A_1 \subseteq A_2
\subseteq ... \subseteq \A$, $B_1 \subseteq B_2 \subseteq ...
\subseteq \B$ such that $\A_0 = \bigcup_{n=0}^{\infty}A_n$, $\B_0
= \bigcup_{n=0}^{\infty}B_n$ form dense unital self-adjoint subrings 
of $\A$ and $\B$, respectively, and pick an increasing
sequence $K_1 \subseteq K_2 \subseteq ... \subseteq G$ of compact
subsets such that $\bigcup_{n=0}^{\infty}K_n = G$.

For each $n$, pick a map $\varphi_n:\B\arrow \Ainfc$ which
satisfies the conditions listed in the statement, with respect to
$A_n,B_n,K_n$ and $\eps = 1/n$.

For each $\varphi_n$, pick a linear self-adjoint lifting (possibly
unbounded)
$$
\tilde{\varphi}_n = (\varphi_n(1),\varphi_n(2),...) : \B \arrow
\lNAc. $$ We claim that we can find a sequence $p_k$ of natural
numbers such that the following hold for all $k \in \N$:

\begin{enumerate}
\item $\| \varphi_{k}(p_k)(b) - \alpha_g(\varphi_{k}(p_k)(b))\| < 1/k$ for all $b\in B_k$, $g \in K_k$.
\item $\|\varphi_{k}(p_k)(1) - 1\| < 1/k$.
\item $\| \varphi_{k}(p_k)(b_1b_2) - \varphi_{k}(p_k)(b_1)\varphi_{k}(p_k)(b_2)\| < 1/k$ for all $b_1,b_2\in B_k$.
\item $\| [\varphi_{k}(p_k)(b),a] \| < 1/k$ for all  $b\in B_k$, $a \in A_k$.
\item $\left | \| \varphi_{k}(p_k)(b) \| - \|b\| \right | < 1/k$ for all $b\in B_k$.
\end{enumerate}
For condition (1), we view $\lNAc$ as the subspace of constant functions in
$\Ch(K_k,\lNAc)$.\\
Consider the linear self-adjoint map
$\tilde{\varphi}'_k: \B \arrow \Ch(K_k,\lNAc)$ given by
$$
\tilde{\varphi}'_k(b)(g) = \bar{\alpha}_g(\tilde{\varphi}_k(b)).
$$
By condition (i) in the statement, we know that
$\|\pi_{K_k}(\tilde{\varphi}_k(b) - \tilde{\varphi}'_k(b))\|<1/k$
for all $b \in B_k$, where $\pi_{K_k}$ denotes the quotient map
onto $\Ch(K_k,\Ainfc)$. Thus for all but finitely many $p \in \N$,
we have that $\| \varphi_{k}(p)(b) - \alpha_g(\varphi_{k}(p)(b))\|
< 1/k$ for all $b \in B_k$.

As for conditions (2) -- (5), we know, from conditions (ii) -- (v)
in the statement, that for all but finitely many $p \in \N$,
conditions (2) -- (5) hold with $p$ standing for $p_k$. Thus, we
can pick some $p_k$ such that all the conditions hold.

Denoting $\pi:\ell^{\infty}(\N,\A) \arrow \Ainf$ the quotient map,
we see that
$$
\pi \circ (\varphi_{1}(p_1),\varphi_{2}(p_2),...)
$$
gives us an isometric unital self-adjoint embedding of $\B_0$ into
$\Ainf$ -- it is a *-ring homomorphism, since all the maps
involved are *-linear, and condition (3) ensures multiplicativity;
it is unital by condition (2), and is isometric by condition (5).
This embedding thus extends to an embedding of $\B$. By conditions
(1) and (4), the image of $\B_0$ commutes with $\A_0$, and is
fixed under the action of $G$. Since $\B_0$ is dense in $\B$ and
$\A_0$ is dense in $\A$, this embedding has the required
properties.

If we do not assume condition (v) of the lemma to hold, the
$\varphi_{k}(p_{k})$ will not necessarily satisfy condition (5)
above and  $\pi \circ (\varphi_{1}(p_1),\varphi_{2}(p_2),...)$
will only extend to a unital $*$-homomorphism. However, if $\B$ is
simple, this $*$-homomorphism will be an embedding since it is
unital.
\end{proof}

\section{Compact Rokhlin actions}

We first recall the definition of the Rokhlin property for finite
groups (see \cite{izumi}).

\begin{Def}  Let $\A$ be a unital separable $C^*$-algebra, $G$ a finite group, and $\alpha:G \arrow \aut(\A)$ an action. $\alpha$ is said to have the \emph{Rokhlin property} if there is a partition of $1_{\A_{\infty}}$ into projections $\{e_g \; | \; g \in G\} \subseteq \Acent$ such that $\bar{\alpha}_g(e_h) = e_{gh}$ for all $g,h \in G$.
\end{Def}

This definition can be generalized in a straightforward way to the
case of compact Hausdorff second-countable groups, as follows.

\begin{Def}  Let $\A$ be a unital separable $C^*$-algebra, $G$ a compact Hausdorff second-countable group,
and $\alpha:G \arrow \aut(\A)$ a strongly continuous action.
$\alpha$ is said to have the \emph{Rokhlin property} if there is a
unital embedding $\Ch(G) \arrow \Acentc$ such that for any $f\in
\Ch(G) \subseteq \Acentc$ (under this embedding), we have
$\bar{\alpha}_g(f)(h) = f(g^{-1}h)$ for all $g,h \in G$.
\end{Def}

\begin{Thm}
\label{compact group main thm}
 Let $\A$ be a separable unital $C^*$-algebra, let $G$ be a compact Hausdorff second-countable group, and
 let $\alpha: G \arrow \aut(\A)$ be an action satisfying the
Rokhlin property. If $\B$ is a unital separable $C^{*}$-algebra which
admits a central sequence of unital embeddings into $\A$, then
$\B$ admits a unital embedding into the fixed point subalgebra of
$\Acent$.
\end{Thm}
\begin{proof}
Fix a unital embedding $\Ch(G) \arrow \Acentc$, as in the
definition of the Rokhlin property (to lighten notation, we shall
view $\Ch(G)$ as embedded in $\Acentc$ in this fashion). Let $B_0
\subseteq \B$, $A_0 \subseteq \A$ be compact subsets.

We shall find a unital embedding $\varphi$ of $\B$ into the fixed point subalgebra of $\Ainfc$ such
that $\|[\varphi(b),a]\|<\eps$ for all $b \in B_0$, $a\in A_0$.
The theorem will follow, then, from Lemma \ref{approximately fixed
commuting embedding}.

We may assume, without loss of generality, that $\alpha_g(A_0)
\subseteq A_0$ for all $g \in G$ (otherwise replace $A_0$ by
$\bigcup_{g \in G} \alpha_g(A_0)$, which is compact, as $A_0$ and
$G$ are compact).

Pick an embedding $\psi:\B \arrow \A$ such that
$\|[\psi(b),a]\|<\eps$ for all $b \in B_0$, $a \in A_0$. Note that
$\|[\alpha_g(\psi(b)),a]\|<\eps$ for all $g\in G$ as well. We view
$\A$ as embedded in $\Ainfc$ in the usual way, and think of $\psi$
as an embedding of $\B$ into $\A \subseteq C^*(\A \cup \Ch(G))
\subseteq \Ainfc$. Note that $C^*(\A \cup \Ch(G)) \cong C(G,\A)$.
We now define $\varphi: \B \arrow C(G,\A)$ by
$$
\varphi(b)(g) = \alpha_g(\psi(b))
$$

It is straightforward to check now that
$\bar{\alpha}_g(\varphi(b)) = \varphi(b)$ for all $g \in G$, and
that $\|[\varphi(b),a]\|<\eps$ for all $b \in B_0$, $a\in A_0$,
giving us the required embedding.
\end{proof}

Recall that a separable unital $C^*$-algebra $\A$ is said to be
\emph{approximately divisible} if there is a finite-dimensional
$C^{*}$-algebra $\B$ with no abelian summands, which admits a unital
embedding into $\MAcent$, or, equivalently, if there is a central
sequence of unital embeddings of $\B$ into $\M(\A)$. We have the
following straightforward consequence of Theorem \ref{compact
group main thm} and Lemma \ref{invariant embedding implies
absorption}.

\begin{Cor} Let $\A$ be a separable unital $C^*$-algebra, let $G$ be a compact Hausdorff second-countable group, and let $\alpha: G \arrow \aut(\A)$ be a strongly continuous action satisfying the Rokhlin property.
\begin{enumerate}
\item If $\D$ is a strongly self-absorbing $C^*$-algebra and $\A$ is $\D$-absorbing, then $\A \times_{\alpha} G$ is $\D$-absorbing.
\item If $\A$ is approximately divisible, then $\A \times_{\alpha} G$ is approximately divisible.
\end{enumerate}
\end{Cor}

\section{Rokhlin Actions of $\Z$}

We recall the definition of the Rokhlin property of an
automorphism (see for instance \cite{kishimoto}).

\begin{Def} Let $\A$ be a unital separable $C^*$-algebra, and $\alpha \in \aut(\A)$. $\alpha$ is said to have the \emph{Rokhlin property} if for any $n$ there is a partition of $1_{\A_{\infty}}$ into projections $e_0,...,e_{n-1},f_0,...,f_n \in \Acent$ such that $\bar{\alpha}(e_j) = e_{j+1}$, $\bar{\alpha}(f_j) = f_{j+1}$ for all $j$, with the convention $e_n = e_0$, $f_{n+1} = f_0$.
\end{Def}

We have the following slight enhancement of the Rokhlin property.
The proof is a simple Cantor-diagonalization type trick, which we
leave to the reader.
\begin{Lemma}
\label{Cantor-diag single auto}
If $\mathcal{C}$ is a separable subspace of $\Acent$, and $\alpha$
is an automorphism of $\A$ satisfying the Rokhlin property,
then the projections as in the Rokhlin property can be chosen to
furthermore commute with $\mathcal{C}$.
\end{Lemma}

In assertion (2) of Theorem \ref{main result} we assumed $\D$ to be $K_{1}$-injective.
In fact, a little less will do:

\begin{Def} We say that a separable unital $C^*$-algebra $\D$ has
\emph{$\oinn$ half-flip} if there is a sequence of unitaries
$u_1,u_2,... \in \U_0(\D \otimes \D)$ such that for all $d \in
\D$, we have $u_n(d\otimes 1)u_n^* \arrow 1 \otimes d$.
\end{Def}

It was shown in \cite{toms-winter}, Proposition 1.13, that if $\D$
is a strongly self-absorbing $C^*$-algebra which is
$K_1$-injective, then it has $\oinn$ half-flip (a review of the
proof shows that, in fact, it is enough to know that $u \otimes u^*
\in \U_0(\D \otimes \D)$ for all $u \in \U(\D)$).\\

The purpose of this section is to prove the following.

\begin{Thm}
\label{single auto main thm} Let $\D$ be a strongly self-absorbing
$C^*$-algebra with an $\oinn$ half-flip, and let $\A$ be a
separable unital $\D$-absorbing algebra. Let $\alpha$ be an
automorphism of $\A$. If $\alpha$ has the Rokhlin property then
$\D$ admits a unital $*$-homomorphism into the fixed point
subalgebra of $\Acent$. In particular, $\A \times_{\alpha}\Z$ is
$\D$-absorbing.
\end{Thm}
Before proving Theorem \ref{single auto main thm}, we need a technical lemma.

\begin{Lemma}
\label{commuting embedding}
Let $\A$, $\B$ be  unital separable $C^*$-algebras. If $\B$ admits
a unital embedding into $\Acent$, and $\mathcal{C}$ is a separable
subalgebra of $\Acent$, then there is a unital embedding of $\B$
into $\Acent \cap \mathcal{C}'$.
\end{Lemma}
\begin{proof}
We denote by $\pi$ the projection $\ell^{\infty}(\N,\A) \arrow
\A_{\infty}$. Let $c_1,c_2,...$ be a dense sequence in
$\mathcal{C}$. We wish to find an embedding of $\B$ into $\Acent$
which commutes with all the elements of this sequence.  Let
$\varphi: \B \arrow \Acent$ be a unital embedding, and let
$\tilde{\varphi} = (\varphi_1,\varphi_2,...):\B \arrow
\ell^{\infty}(\N,\A)$ be a linear lifting. Note that for any
subsequence $(\varphi_{k_1},\varphi_{k_2},...)$, composition with
$\pi$ will give us another unital embedding of $\B$ into $\Acent$.
We can thus use a Cantor-diagonalization type argument to
construct such a subsequence, such that the image $\pi \circ
(\varphi_{k_1},\varphi_{k_2},...)$ will commute with
$\mathcal{C}$. More specifically, note first that for any $a \in
\A$, $b \in \B$, we have that $ \|[\varphi_n(b),a]\| \arrow 0$.
For each $c_k$ we pick a lifting $\tilde{c_k} =
(c_k(1),c_k(2),...) \in \ell^{\infty}(\N,\A)$. Let $b_1,b_2,...$
be a dense sequence in $\B$. For each $n$, choose $k_n$ such that
$\|[\varphi_{k_n}(b_j),c_i(n)]\|<1/n$ for all $i,j \leq n$ (and we
may assume that this sequence $k_n$ is increasing). Thus,  $\pi
\circ (\varphi_{k_1},\varphi_{k_2},...)$ gives us an embedding
with the required properties.
\end{proof}

We can now prove Theorem \ref{single auto main thm}.

\begin{proof}[Proof of Theorem \ref{single auto main thm}]
It suffices, by Lemmas \ref{approximately fixed commuting
embedding} and \ref{invariant embedding implies absorption}, to
prove that for any finite subset $F \subseteq \D$ and $\eps>0$
there is a unital embedding
$$\varphi:\D \arrow \Acent$$ such that
\[
\|\bar{\alpha}(\varphi(x)) - \varphi(x)\| < \eps
\]
for all $x \in F$.
Let us fix such $F,\eps$.

Fix an embedding $\iota:\D \arrow \Acent$. $C^*\left (
\bigcup_{k=-\infty}^{\infty}\bar{\alpha}^k(\iota(D)) \right )$ is
separable, so by Lemma \ref{commuting embedding}, there is a
unital embedding of $\eta:\D \arrow \Acent$ which commutes with
it. Let $$\textstyle \B := C^*\left ( \eta(\D) \cup
\bigcup_{k=-\infty}^{\infty} \bar{\alpha}^k(\iota(D)) \right ) \,
; $$ note that $\B \cong \B \otimes \D$.

Choose a unitary $w \in \U_0(\D \otimes \D)$ such that
$$\|w(x\otimes 1)w^*-1 \otimes x\| < \frac{\eps}{4}$$ for all $x \in F$. $w$
can be connected to $1_{\D \otimes \D}$ via a rectifiable path.
Let $L$ be the length of such a path. Choose $n$ such that
$L\|x\|/n <\eps/8$ for all $x \in F$.

Consider the embeddings $\iota, \bar{\alpha}^n \circ \iota :\D \arrow
\B$. Define embeddings $$\rho,\rho':\D \otimes \D \arrow \B$$ by
$$\rho(x \otimes y) = \iota (x) \eta(y), \;\; \rho'(x\otimes y) =
\bar{\alpha}^n(\iota(x))\eta(y) \, .$$ Pick unitaries $1=w_0,w_1,...,w_n = w \in
\U_0(\D \otimes \D)$ such that $\|w_k-w_{k+1}\| \leq L/n$ for
$k=0,...,n-1$. Now, let $u_k = \rho(w_k)^*\rho'(w_k)$. Note that
$$\|u_k - u_{k+1}\| \leq \frac{2L}{n}$$ for $k=0,...,n-1$, that $u_0 = 1$,
and that $$\|u_n \bar{\alpha}^n(\iota(x)) u_n^* - \iota(x) \| <
\frac{\eps}{2}$$ for all $x\in F$.

Similarly, we choose unitaries $1=v_0 ,v_1,...,v_{n+1} \in \B$
such that $$\|v_k-v_{k+1}\| \leq \frac{2L}{n+1}$$ for $k=0,...,n$
and
$$\|v_{n+1}\bar{\alpha}^{n+1}(\iota(x))v_{n+1}^* - \iota(x)\| < \frac{\eps}{2}$$
for all $x \in F$.

We use the Rokhlin property (and Lemma \ref{Cantor-diag single auto}) to find a partition of $1_{\A_{\infty}}$ into projections
$$
e_i,f_j \in \Acent\cap \B' \cap
\{\bar{\alpha}^{k-n}(u_k),\bar{\alpha}^{l-n-1}(v_j) \, | \, k=1,...,n-1,
l=1,...,n\}' \, ,
$$
$i=0, \ldots,n-1$, $j=0, \ldots,n$, such that $\bar{\alpha}(e_j) =
e_{j+1}$, $\bar{\alpha}(f_j) = f_{j+1}$ for all $j$, with the
convention $e_n = e_0$, $f_{n+1} = f_0$. We then define
$\varphi:\D \arrow \Acent$ by
$$
\varphi(x) = \sum_{k=0}^{n-1}e_k\bar{\alpha}^{k-n}(u_k)\bar{\alpha}^k(\iota(x))\bar{\alpha}^{k-n}(u_k^*)e_k + \sum_{k=0}^{n}f_k\bar{\alpha}^{k-n-1}(v_k)\bar{\alpha}^k(\iota(x))\bar{\alpha}^{k-n-1}(v_k^*)f_k .
$$
By the choice of the projections $e_0,...,e_{n-1},f_0,...,f_n$,
$\varphi$ is indeed a unital homomorphism, and a simple
computation shows that for any $x \in F$, we have that
$$\|\bar{\alpha}(\varphi(x)) - \varphi(x)\| < \frac{\eps}{2} + \frac{4L\|x\|}{n} <
\eps \, ,$$ as required.
\end{proof}

\section{Rokhlin flows}

Recall the definition of a Rokhlin flow from \cite{kishimoto flows}:
\begin{Def}
A strongly continuous action $\alpha: \R \to \aut(\A)$ on a unital
separable $C^{*}$-algebra $\A$ is said to have the Rokhlin
property if, for any $p \in \R$, there is a unitary $v \in
\Acentc$ such that $\bar{\alpha}_{t}(v)= e^{itp} \cdot v$.
\end{Def}

Our main theorem for this section is the following.

\begin{Thm}
\label{main theorem for flows}
Let $\A$ and $\D$ be separable unital $C^{*}$-algebras with $\D$ strongly self-absorbing with an $\oinn$ half-flip, and suppose $\A$ is $\D$-absorbing; let $\alpha: \R \to \aut(\A)$ be a Rokhlin flow. \\
Then, there is a unital $*$-homomorphism
\[
\varphi: \D \to \Acent
\]
which is invariant under $\bar{\alpha}$. In particular, $\A
\times_{\alpha} \R$ is $\D$-absorbing.
\end{Thm}

It was shown in \cite{kishimoto flows} that if $\A$ is simple and
purely infinite with a Rokhlin flow $\alpha$, then the crossed
product is simple and purely infinite again. In the case  where
$\A$ is simple, unital and nuclear, for $\D=\OO_{\infty}$ our
result coincides with Kishimoto's theorem (although the proofs are
very different), since Kirchberg has shown that being purely
infinite and absorbing $\OO_{\infty}$ are equivalent conditions when $\A$ is
nuclear (see \cite{kirchberg book,kirchberg-phillips}).

The basic idea of the proof of Theorem \ref{main theorem for
flows} is quite similar to that of Theorem \ref{single auto main
thm} (we organize it differently, though). The main difference
stems from the fact that the induced action of $\R$ on $\Acent$ is
generally discontinuous -- a problem which is irrelevant in the
case of $\Z$. Since embeddings into $\Acent$ arise as approximate
embeddings into $\A$, we work here with maps into $\A$ and into
$C^*(\A,v) \subseteq \Ainfc$, where $v$ is a unitary arising from
the Rokhlin property, which behave approximately like their
analogues from the previous section. This allows us to circumvent
the continuity problems, but at the cost of having messier
estimates.The algebra $C^*(v)$, which is in fact isomorphic to
$\Ch(\T)$, here will play an analogous role to that of $\Ch(G)$ in
the compact case.

\begin{Lemma}
\label{cutdown lemma} Let $\D$ be a strongly self-absorbing
$C^*$-algebra, and let $\A$ be a unital separable
$C^*$-algebra. For any compact subset $K \subset \A$ and $\eps>0$
there is a unital expectation $\theta:\A \arrow \A$ such that
$\|x-\theta(x)\|<\eps$ for all $x\in K$, and such that $\D$ admits
a unital embedding into $\A \cap \theta(\A)'$.
\end{Lemma}
\begin{proof}
We have that $$ \A \cong \A \otimes \D^{\otimes \infty} =
\lim_{\longrightarrow} (\A\otimes \D^{\otimes k}) \, .
$$
Under this identification, there is some $k$ such that
$\textrm{dist}(K,\A \otimes \D^{\otimes k} \otimes 1 )< \eps$,
where $1$ here denotes the identity of the copy of $\D^{\otimes
\infty}$ sitting in $\A \otimes \D^{\otimes \infty}$ as $1_{\A}
\otimes 1_{\D^{\otimes k}} \otimes \D^{\otimes \infty}$. If we fix
some arbitrary state $\tau$ of $\D^{\otimes \infty}$, the map
$\theta = \id_{\A \otimes \D^{\otimes k}} \otimes \tau$ will
satisfy the required properties.
\end{proof}

In the following lemma, we think of $\Ch(\T,\A)$ as the space of
$\A$-valued $2M$-periodic functions on $\R$ (where $M$ is selected
in the lemma).

\begin{Lemma}
\label{periodization lemma} Let $\A$ and $\D$ be separable unital
$C^{*}$-algebras; suppose that $\D$ is strongly self-absorbing
with an $\oinn$ half-flip and that $\A\cong \A\otimes \D$. Let
$\alpha$ be a strongly continuous action of $\R$ on $\A$. For any
$\eps>0$, any compact sets $D_0 \subseteq \D$, $A_0 \subseteq \A$
and any $\mu>0$ there are an $M>\mu$ and a completely positive
contraction $\beta:\D \arrow \Ch(\T,\A)$, with the following
properties:
\begin{enumerate}
\item \label{PL - alpha periodic} $\|\alpha_t(\beta(d)(s-t)) - \beta(d)(s)\|< \eps$
\item \label{PL - unital} $\|\beta(1) - 1\| < \eps$
\item \label{PL - multiplicative} $\|\beta(dd') - \beta(d)\beta(d')\|<\eps$
\item \label{PL - commuting} $\|[\beta(d),a]\|<\eps$ (where $a$ is thought of as a constant
function on $\T$)
\end{enumerate}
for all $d,d' \in D_0$, $a \in A_0$, $t \in [-\mu,\mu]$, $s \in
[-M,M]$.
\end{Lemma}

As the proof of this lemma is somewhat lengthy, we shall first
show how it is used to prove Theorem \ref{main theorem for flows}.

\begin{proof}[Proof of Theorem \ref{main theorem for flows}.]
Fix compact sets $D_0 \subseteq \D$, $A_0 \subseteq \A$,
$[-\mu,\mu] \subseteq \R$ ($\mu >0$) and $\eps>0$. By Lemmas
\ref{invariant embedding implies absorption} and
\ref{approximately fixed commuting embedding}, it will suffice to
find a completely positive contraction $\varphi:\D \arrow \Ainfc$
such that, for any $d,d' \in D_0$, $a \in A_0$, $t \in [-\mu,\mu]$
we have that the expressions $\|\bar{\alpha}_t(\varphi(d)) -
\varphi(d)\|$,   $\|\varphi(1)-1\|$, $\|\varphi(d)\varphi(d') -
\varphi(dd')\|$  and $\|[\varphi(d),a]\|$ are all bounded above by
$\eps$. (Since $\D$ is simple, we do not have to verify condition
(v) of Lemma \ref{approximately fixed commuting embedding}.)

By the Rokhlin property of $\alpha$, for any $M>0$ there is a
unitary $v \in \Acentc$ satisfying
\[
\bar{\alpha}_{t}(v) = e^{-\pi i t/M}\cdot v \; \forall \, t \in \R
\, .
\]
We may define a $*$-homomorphism
\[
\sigma: \Ch(\T) \otimes \A  \to C^{*}(v, \A) \subseteq \Ainfc
\]
by
\[
\sigma(f \otimes x) = f(v)x \;\;\; \forall \, x \in \A, f \in
\Ch(\T) \, .
\]
We then have
\[
\bar{\alpha}_{t} \circ \sigma(f \otimes x) = \sigma(f_{- t}
\otimes \alpha_{t}(x))
\]
for all $t \in \R$, $f \in \Ch(\T)$ and $x \in \A$, where $f_{-
t}$ denotes the function $f$ translated by $-t$. We write $\tilde{\alpha}_{t}$ for the action on $\Ch(\T) \otimes A$ given by $f \otimes x \mapsto f_{-t} \otimes \alpha_{t}(x)$.  As before, we identify
$\Ch(\T) \otimes A \cong \Ch(\T,\A)$ with the space of $\A$-valued $2M$-periodic functions
on $\R$.

Select $M$, $\beta$
as in Lemma \ref{periodization lemma} and define $$ \varphi =
\sigma \circ \beta \ .  $$ We only check that
$\|\bar{\alpha}_t(\varphi(d)) - \varphi(d)\| < \eps$ for $d \in
D_{0}$, $t \in [-\mu,\mu]$; the other requirements above are
straightforward to verify:
\begin{eqnarray*}
\|\bar{\alpha}_t(\varphi(d)) - \varphi(d)\| & \le & \|\tilde{\alpha}_{t} \circ \beta(d) - \beta(d)\| \\
& = & \sup_{s \in [-M,M]} \|(\tilde{\alpha}_{t} \circ \beta(d))(s) - \beta(d)(s)\| \\
& = & \sup_{s \in [-M,M]} \|\alpha_{t}(\beta(d)(s-t)) - \beta(d)(s)\| \\
& \stackrel{\ref{periodization lemma}(1)}{<} & \eps \ .
\end{eqnarray*}
\end{proof}

It remains to prove Lemma \ref{periodization lemma}.
\begin{proof}[Proof of Lemma \ref{periodization lemma}.]
We may assume, without loss of generality, that $\|x\|\leq 1$ for
all $x\in D_0 \cup A_0$, and that $1_{\D} \in D_0$, $1_{\A} \in A_0$ (in
particular, we have $D_0 \subseteq D_0^2 = \{xy \; | \; x,y \in
D_0\})$.

Since $\D$ is strongly self-absorbing with $\oinn$ half-flip,
there is a unitary $w \in \U_0(\D \otimes \D)$ such that
\begin{equation} \label{eq a}
\|w(d \otimes 1_{\D})w^{*} - 1_{\D} \otimes d\| <
\frac{\varepsilon}{40} \;\; \forall \, d \in D_0^2 \,  \, .
\end{equation}
There is thus a continuous path of some finite length $L$ in the unitary group of $\D \otimes \D$ connecting $w$ and $1_{\D \otimes \D}$.\\
Fix some $M > \max\{\mu,40\mu L/\eps\}$. Pick a continuous path
$(w_t)_{t\in [-M,M]} \subseteq \U_0(\D \otimes \D)$,  such that
$w_{-M} = 1_{\D \otimes \D}$, $w_M = w$ and
\[ \|w_{t}-w_{s}\|
\le \frac{L|t-s|}{2M}  \; \; \forall \, s,t \in [-M,M] \, .
\]

We find $x_1,...,x_p,y_1,...,y_p \in \D$ of norm $\leq 1$, and
$\lambda_j(i)\in\C$, $j = 1,...,p$, $i=1,...,P$ for some $P$, such
that for any $t \in [-M,M]$ there is an $i \in \{1, \ldots, P\}$ such that
\begin{equation} \label{eq b}
 \|\sum_{j=1}^p\lambda_j(i)x_j \otimes y_j - w_t\| <
\frac{\eps}{400} \end{equation}
and $\|\sum_{j=1}^p\lambda_j(i)x_j
\otimes y_j\|\leq 1$ for all $i$. Let $\Lambda =
\max\{|\lambda_j(i)| +1 \; | \; j\leq p, i\leq P\}$.

Pick a unital embedding $\psi:\D \arrow \A$ such that
\begin{equation} \label{eq c}
\|[\psi(x),a]\|<\frac{\eps}{30p\Lambda}<\frac{\eps}{30}
\end{equation}
for all $a \in \bigcup_{t=-2M}^{2M}\alpha_t(A_0)$ and all $x \in
\{x_1,...,x_p,x_1^*,...,x_p^*\} \cup D_0$. Use Lemma \ref{cutdown
lemma} to find an expectation map $\theta: \A \arrow \A$ such that
\begin{equation} \label{eq d}
\|\theta(x) - x\|<\frac{\eps}{400p^2\Lambda^2} < \frac{\eps}{400}
\end{equation}
for all
$$
x \in\left ( \bigcup_{t=-3M}^{3M}\alpha_t(\psi(
\{x_1,...,x_p,x_1^*,...,x_p^*\} \cup D_0) \cup A_0) \right )^3 \,
,
$$
and a unital embedding $\eta: \D \arrow \A \cap\theta(\A)'$. \\
Notice that for all $x,y,z \in
\bigcup_{t=-M}^M\alpha_t(\psi( \{x_1,...,x_p,x_1^*,...,x_p^*\}
\cup D_0) \cup A_0)$, we have
\begin{equation} \label{eq e}
\|\theta(xyz) - \theta(x)\theta(y)\theta(z)\| <
\frac{4\eps}{400p^2\Lambda^2} < \frac{4\eps}{400} \, .
\end{equation}
Define unital completely positive maps $ \varrho,
\varrho': \D \otimes \D \to \A $ by
\[
\varrho(d_{1}\otimes d_{2}):=  \theta(\alpha_{-M}(\psi (d_{1})))
\cdot \eta(d_{2}) \;\; , \;\; \varrho'(d_{1}\otimes d_{2}):=
\theta(\alpha_{M}(\psi(d_{1}))) \cdot \eta(d_{2}) \] for
$d_{1},d_{2} \in \D$.  For any $t \in [-M,M]$ and $d_1,d_2 \in
D_0^2$, we find $i\in \{1, \ldots , P\}$ such that $\|\sum_{j=1}^p\lambda_j(i)x_j
\otimes y_j - w_t\| < \eps/400$; we see then that
\begin{eqnarray} 
\label{eq f}
\lefteqn{\|\varrho(w_t)\varrho(d_1 \otimes d_2) \varrho(w_t^*) -
\varrho(w_t (d_{1} \otimes d_2) w_t^*)\|}
\\\nonumber
& \underbrace{<}_{(\ref{eq b})} & \left \| \varrho \left (
\sum_{j=1}^p\lambda_j(i)x_j \otimes y_j \right ) \varrho(d_1
\otimes d_2) \varrho\left (
\sum_{j=1}^p\lambda_j(i) x_j \otimes y_j \right ) ^*  \right .
\\\nonumber
 & & - \left . \varrho \left ( \left (
\sum_{j=1}^p\lambda_j(i)x_j \otimes y_j \right ) (d_1 \otimes d_2)
\left (\sum_{j=1}^p\lambda_j(i)x_j \otimes y_j\right )^*
\right ) \right \|  + \frac{4\eps}{400}
\\\nonumber
& \leq & \sum_{j,k=1}^p|\lambda_j(i)\overline{\lambda_k(i)}|  \| (
\theta(\alpha_{-M}(\psi(x_j)))\theta(\alpha_{-M}(\psi(d_1)))\theta(\alpha_{-M}(\psi(x_k^*)))
\\\nonumber
& &
- \theta(\alpha_{-M}(\psi(x_j))\alpha_{-M}(\psi(d_1))\alpha_{-M}(\psi(x_k^*)))
)\|\| \eta(y_jd_2y_k^*)  \| + \frac{4\eps}{400}
\\\nonumber
& \underbrace{<}_{(\ref{eq e})} & \frac{4\eps}{400} + \frac{4\eps
p^2\Lambda^2}{400p^2\Lambda^2}
\\ \nonumber
& = & \frac{\eps}{50}
\end{eqnarray}
and the same estimate holds if we replace $\varrho$ by $\varrho'$,
and if we interchange $w_t$ and $w_t^*$. In particular, setting
$d_1 = d_2 = 1$, we see that for all $t\in [-M,M]$
\begin{equation} \label{eq g}
\|\varrho(w_t^*)\varrho(w_t) - 1\|, \|\varrho(w_t)\varrho(w_t^*) -
1\|, \|\varrho'(w_t^*)\varrho'(w_t) - 1\|,
\|\varrho'(w_t)\varrho'(w_t^*) - 1\| < \frac{\eps}{50} \, .
\end{equation}

Similarly, we estimate for any $s \in [-M,M]$, $t \in [-2M,2M]$
and $a \in A_0$ (for a suitable $i \in \{1, \ldots, P\}$):
\begin{eqnarray} \label{eq h}
\|[\varrho(w_s),\alpha_t(a)]\|
& \underbrace{<}_{(\ref{eq b})} &
\left \| \left[ \varrho \left ( \sum_{j=1}^p\lambda_j(i)x_j
\otimes y_j \right ) , \alpha_t(a)
\right ] \right \| + \frac{2\eps}{400}
\\\nonumber
& \underbrace{<}_{(\ref{eq d})} &
\left \| \left[ \varrho \left (
\sum_{j=1}^p\lambda_j(i)x_j \otimes y_j \right ) ,
\theta(\alpha_t(a)) \right ] \right \| + \frac{4\eps}{400}
\\\nonumber
& \leq &  p\Lambda\max_{j\leq
p}\|[\theta(\alpha_{-M}(\psi (x_j))),\theta(\alpha_t(a))]\| + \frac{4\eps}{400}
\\\nonumber
& \underbrace{\leq}_{(\ref{eq e})} & p\Lambda\max_{j\leq
p}\|\theta([\psi (x_j),\alpha_{t+M}(a)])\| + \frac{12\eps}{400}
\\\nonumber
& \underbrace{<}_{(\ref{eq c})} &  \frac{\eps}{30} + \frac{12\eps}{400}
\\\nonumber
& < & \frac{\eps}{15}
\end{eqnarray}
and the same estimate holds with $w_s^*$ instead of $w_s$, and
with $\varrho'$ instead $\varrho$.

Define a continuous path $(u_{t})_{t\in[-M,M]} \subseteq \A$ by
\[
u_{t}:= \varrho(w_{t})^{*}\varrho'(w_{t}) \, .
\]
The $u_t$'s are not necessarily unitary, but one checks that they
satisfy
\begin{equation} \label{eq i}
\|u_tu_t^* - 1\|, \|u_t^*u_t - 1\| <
\frac{2\eps}{50} = \frac{\eps}{25}\, .
\end{equation}
Note also that
\begin{equation} \label{eq j}
\|u_t - u_s\| \le \frac{2L}{2M} \cdot |t-s| \; \forall \, s,t \in
[-M,M]
\end{equation}
and that $u_{-M} = 1$. Let $d \in D_0^2$. Using that
$$\|\alpha_{M}(\psi(d)) - \varrho'(d \otimes 1_{\D})\| \ , \
\|\alpha_{-M}(\psi(d)) - \varrho(d \otimes 1_{\D})\|<\eps/400 \ , $$
we check that
\begin{eqnarray}\label{eq k}
\lefteqn{ \|u_{M}\alpha_{M}(\psi(d)) u_{M}^{*} -
\alpha_{-M}(\psi(d))\|}
\\\nonumber
& <  &  \|\varrho(w_{M})^{*}\varrho'(w_{M}) \varrho'(d \otimes
1_{\D}) \varrho'(w_{M})^{*}\varrho(w_{M}) - \varrho(d \otimes
1_{\D})\| + \frac{2\eps}{400}
\\\nonumber
& \underbrace{<}_{(\ref{eq f})} &
\|\varrho(w_{M})^{*}\varrho'(w_{M}(d \otimes 1_{\D})
w_{M}^*)\varrho(w_{M}) - \varrho(d \otimes 1_{\D})\| +
\frac{2\eps}{400} + \frac{\eps}{50}
\\\nonumber
& \underbrace{<}_{(\ref{eq a})} & \|\varrho(w_{M})^{*}
\varrho(1_{\D} \otimes d) \varrho(w_{M}) - \varrho(d \otimes
1_{\D})\| + \frac{2\eps}{400} + \frac{\eps}{50} + \frac{\eps}{40}
\\\nonumber
&<& \frac{2\eps}{400} + 2\cdot\left (\frac{\eps}{50} +
\frac{\eps}{40}\right)
\\\nonumber
& < & \frac{\eps}{10} \ .
\end{eqnarray}
Furthermore,
\begin{eqnarray} \label{eq l}
\|[\alpha_{t}(a),u_{s}]\|
& = &  \|[\alpha_{t}(a), \varrho(w_{s}^{*})\varrho'(w_{s})]\|
\\\nonumber
& \le &  \|[\alpha_{t}(a), \varrho(w_{s})]\| + \|[\alpha_{t}(a), \varrho'(w_{s})]\|
\\\nonumber
 & \underbrace{<}_{(\ref{eq h})} &   \frac{2\varepsilon}{15}
\end{eqnarray}
for $a \in A_0$, $s \in [-M,M]$, $t \in [-2M,2M]$, and the same estimate holds for $u_{s}^{*}$ in place of $u_{s}$.

For each $d \in \D$, $t\in [-M,M]$, we define
\[
h(d,t):= \alpha_{t-M}(u_{t}) \alpha_{t}(\psi(d))
\alpha_{t-M}(u_{t}^{*}) \, ,
\]
then
\[
h(d,M) = u_{M} \alpha_{M}(\psi(d))u_{M}^{*}
\]
and
\[
h(d, -M) = \alpha_{-M}(\psi(d)) \; \forall \, d \in \D \, ,
\]
whence
\begin{equation} \label{eq lb}
\|h(d,M) - h(d, -M)\| < \frac{\varepsilon}{10} \; \forall \, d \in
D_0^2 \, .
\end{equation}
Note that for any $d_1,d_2$ in the unit ball of $\D$ and for any
$t \in [-M,M]$ we have that
\begin{equation} \label{eq m}
\|h(d_1,t)h(d_2,t) - h(d_1d_2,t)\| \leq \|u_t^*u_t - 1\|
\underbrace{<}_{(\ref{eq i})} \frac{\eps}{25} \, .
\end{equation}
For any $a \in A_0$, $d \in D_0$, $t \in [-M,M]$, we have
\begin{eqnarray} \label{eq n}
\|[h(d,t),a]\|& \leq &
\|[\alpha_{t-M}(u_{t}),a]\|+\|[\alpha_{t}(\psi(d)),a]\|+\|[\alpha_{t-M}(u_{t}^{*}),a]\|
\\\nonumber
& = &
\|[u_{t},\alpha_{M-t}(a)]\|+\|[\psi(d),\alpha_{-t}(a)]\|+\|[u_{t}^{*},\alpha_{M-t}(a)]\|
\\\nonumber
 & \underbrace{<}_{(\ref{eq c}),(\ref{eq l})} &
 \frac{2\eps}{15}+\frac{\eps}{30}+\frac{2\eps}{15}
\\\nonumber
& < & \frac{\eps}{3} \, .
\end{eqnarray}

We can view $h$ as a completely positive contraction from $\D$ to
$\Ch([-M,M],\A)$. $\Ch(\T,\A)$ may be identified with the subalgebra of
$\Ch([-M,M],\A)$ consisting of functions which agree on $-M$ and
$M$. We shall now perturb $h$ so as to ensure that its range is in
$\Ch(\T,\A)$. For each $\delta>0$, define a continuous function
$g_{\delta} \in \Ch([-M,M])$ by
\[
g_{\delta}(t):= \left\{
\begin{array}{ll}
0, & t = -M \\
1, & -M+\delta \le t \le M \\
\mbox{linear } & \mbox{elsewhere} \, ;
\end{array}
\right.
\]
Now, let
\[
h_{\delta}(d,t):= g_{\delta}(t) \cdot h(d,t) + (1 - g_{\delta}(t))
\cdot h(d,M) \, .
\]
For any $d \in \D$, $h_{\delta}(d,t)$ is a continuous function of
$t$, and satisfies $h_{\delta}(d,-M) = h_{\delta}(d,M)$.
$h_{\delta}$, thus, can be regarded as a map $ h_{\delta}:\D \arrow
\Ch(\T,\A)$. This map can readily be seen to be a completely positive
contraction.

Fix $\delta>0$ such that
\[
\|h(d,t) - h_{\delta}(d,t)\| < \frac{\varepsilon}{4} \;\; \forall
\, d \in D_0^2 , \, t \in [-M,M] \, ;
\]
this is possible by (\ref{eq lb}).
We let $$\beta := h_{\delta} \, .$$
The proof will be complete once we
check that $\beta$ satisfies the four conditions required in the
statement. \\ \underline{Condition \ref{PL - alpha periodic}}: We
wish to show that
\[
\|\alpha_{t}(\beta(d,s- t)) - \beta(d,s)\| < \eps \;\; \forall \,
d \in D_0, \, s \in [-M,M], \, t \in [-\mu,\mu] \, .
\]
It will clearly be enough to consider $t \in [0,\mu]$. (The case $t<0$ will follow with $t$ replaced by $-t$ and $s$ replaced by $s-t$.) Fix $d\in
D_0$.
\begin{itemize}
\item[\emph{Case I:}] $s \in [-M+ t,M]$. We have
\begin{eqnarray*}
\lefteqn{\|\alpha_{t}(h(d,s - t)) - h(d,s)\| }\\
& = & \|\alpha_{s - M}(u_{s - t}) \alpha_{s}(\psi(d)) \alpha_{s - M}(u_{s -t}^{*}) -
\alpha_{s - M}(u_{s}) \alpha_{s}(\psi(d)) \alpha_{s - M}(u_{s}^{*})\| \\
& \leq & 2\|u_{s-t} - u_s\|
\underbrace{\leq}_{(\ref{eq j})}
\frac{2 L|t|}{M} \leq \frac{2L\mu}{M} < \frac{\eps}{20}\, ,
\end{eqnarray*}
from which follows that
\[
\|\alpha_{t}(\beta(d,s-t)) - \beta(d,s)\| <
2\cdot\frac{\eps}{4}+\frac{\eps}{20} < \eps \, .
\]
\item[\emph{Case II:}] $s \in [-M, -M+ t]$. In this case, there are $0 \le t_{0}, t_{1}\le t$ such
that $t=t_{0}+t_{1}$ and $s = -M + t_{0}$.  Note that
$\|\beta(d,s-t) - h(d,M -  t_{1}) \| < \eps/4$. Thus
\begin{eqnarray*}
\lefteqn{\|\alpha_{t}(\beta(d,s-t)) - \beta(d,s)\|}\\
& \leq & \|\alpha_{t}(h(d,M -  t_{1})) - h(d,s)\| + \frac{2\eps}{4}\\
& \leq & \| \alpha_{t_{0}}(h(d,M)) - h(d,s)\| +
\|\alpha_{t_1}(h(d,M -  t_{1})) - h(d,M)\|
+ \frac{2\eps}{4}  \\
& < & \| \alpha_{-M+t_{0}}(\psi(d)) - h(d,s)\| + \frac{3
\varepsilon}{4} + \frac{\eps}{20}\\
& \leq &  2\|u_{-M} - u_{s}\| +
\frac{3 \varepsilon}{4} + \frac{\eps}{20} \\ &<&  \frac{\eps}{20}
+ \frac{3 \varepsilon}{4} + \frac{\eps}{20} < \eps \, .
\end{eqnarray*}
\end{itemize}
\underline{Condition \ref{PL - unital}}: Indeed,
\[
\|\beta(1) - 1\| <  \sup_{t \in [-M,M]}\|h(1,t) - 1\| +
\frac{\eps}{4} =  \sup_{t \in [-M,M]}\|\alpha_{t-M}(u_tu_t^*) -
1\| + \frac{\eps}{4} < \frac{\eps}{25}+\frac{\eps}{4} < \eps \ .
\]
\underline{Condition \ref{PL - multiplicative}}: Let $d_1,d_2 \in
D_0$. We have
\begin{eqnarray*}
\|\beta(d_1)\beta(d_2) - \beta(d_1d_2)\|
& \le & \sup_{t\in[-M,M]} \|\beta(d_1,t)\beta(d_2,t) - \beta(d_1d_2,t))\| \\
& < & \sup_{t\in[-M,M]}  \|h(d_1,t)h(d_2,t) - h(d_1d_2,t))\| +
\frac{3\varepsilon}{4} \\
& \underbrace{<}_{(\ref{eq m})}  & \frac{\varepsilon}{25} +
\frac{3\varepsilon}{4} \\
& <  & \eps \ .
\end{eqnarray*}
\underline{Condition \ref{PL - commuting}}: For each $d \in D_0$,
$a\in A_0$, we have
\[
\|[\beta(d),a]\| < \sup_{t\in [-M,M]}\|[h(d,t),a]\| +
\frac{2\eps}{4}
 \underbrace{<}_{(\ref{eq n})}
 \frac{\eps}{3} + \frac{2\eps}{4} < \eps \ .
\]

\end{proof}

\end{document}